\documentclass[12pt,thmsa]{article}
\usepackage{amssymb}
\usepackage{sw20bams}

\input tcilatex
\begin{document}

\author{Steven R. Finch}
\title{CLT\ Variance Associated with Baxendale's SDE}
\date{September 30, 2008}
\maketitle

\begin{abstract}
Simple analysis of the leftmost eigenvalue of Ince's equation (a boundary
value problem with periodicity) resolves an open issue surrounding a
stochastic Lyapunov exponent. Numerical verification is also provided.
\end{abstract}

\footnotetext{
Copyright \copyright\ 2008 by Steven R. Finch. All rights reserved.}Let $a>b$
and $\sigma >0$. Consider the stochastic differential equation (SDE) 
\[
\begin{array}{ccc}
dX_t=\left( 
\begin{array}{cc}
a & 0 \\ 
0 & b
\end{array}
\right) X_t\,dt+\sigma \left( 
\begin{array}{cc}
0 & -1 \\ 
1 & 0
\end{array}
\right) X_t\circ dW_t, &  & X_t\in \Bbb{R}^2
\end{array}
\]
where $W_t$ is scalar Brownian motion with unit variance and $\circ $
denotes the use of Stratonovich calculus. A measure of the stability of such
a system is provided by the (almost-sure) Lyapunov exponent 
\[
\ell =\lim_{t\rightarrow \infty }\frac 1t\ln |X_t| 
\]
for $X_0\neq 0$. Baxendale \cite{Bx, KP} computed a formula for $\ell $ and,
further, proved an associated central limit theorem (CLT) 
\[
\lim_{t\rightarrow \infty }\limfunc{P}\left( \frac{\dfrac 1t\ln |X_t|-\ell }{%
\dfrac s{\sqrt{t}}}\leq v\right) =\frac 1{\sqrt{2\pi }}\dint\limits_{-\infty
}^v\exp \left( -\frac{u^2}2\right) du. 
\]
No formula has been known for the variance $s^2$ until now. Our calculation
of $s^2$ is based on the boundary value problem (BVP) 
\[
y^{\prime \prime }(x)+c\sin (2x)y^{\prime }(x)+\left( \lambda -\mu \,c\cos
(2x)\right) y(x)=0, 
\]
\[
y^{\prime }(0)=y^{\prime }(\pi )=0 
\]
due to Ince \cite{Ic}, where $c=-(a-b)/\sigma ^2$ and $\mu \approx 0$. More
precisely, if $\lambda (\mu )$ is the leftmost eigenvalue of the Ince BVP,
given $\mu $, then 
\[
\begin{array}{ccc}
\ell =\dfrac{a+b}2-\dfrac{\sigma ^2}2\lambda ^{\prime }(0), &  & s^2=-\dfrac{%
\sigma ^2}2\lambda ^{\prime \prime }(0).
\end{array}
\]
In the next three sections, we discuss how Ince's equation arises and the
details of computing $\lambda ^{\prime }(0)$ and $\lambda ^{\prime \prime
}(0)$. The final section is devoted to numerical verification of the
preceding.

\section{Equivalence}

Let $\gamma /\sigma ^2$ denote the rightmost eigenvalue of the differential
operator \cite{Bx} 
\[
\frac 12\frac d{dx^2}-\frac{a-b}{2\sigma ^2}\sin (2x)\frac d{dx}+\frac \mu
{\sigma ^2}(a\cos ^2x+b\sin ^2x) 
\]
(which is obtained, in part, by projecting the solution $X_t$ of Baxendale's
SDE onto the unit circle). By the double angle formula for cosine, the BVP\
we wish to study is 
\[
y^{\prime \prime }(x)-\frac{a-b}{\sigma ^2}\sin (2x)y^{\prime }(x)+\left[ -%
\frac{2\gamma }{\sigma ^2}+\frac \mu {\sigma ^2}\left( (a+b)+(a-b)\cos
(2x)\right) \right] y(x)=0 
\]
hence $c=-(a-b)/\sigma ^2$ follows, as does 
\[
\lambda =-\frac{2\gamma }{\sigma ^2}+\frac{a+b}{\sigma ^2}\mu ; 
\]
hence 
\[
\begin{array}{ccccc}
\lambda ^{\prime }=-\dfrac{2\gamma ^{\prime }}{\sigma ^2}+\dfrac{a+b}{\sigma
^2} &  & \therefore &  & \ell =\dfrac{a+b}2-\dfrac{\sigma ^2}2\lambda
^{\prime }(0),
\end{array}
\]
\[
\begin{array}{ccccc}
\lambda ^{\prime \prime }=-\dfrac{2\gamma ^{\prime \prime }}{\sigma ^2} &  & 
\therefore &  & s^2=-\dfrac{\sigma ^2}2\lambda ^{\prime \prime }(0).
\end{array}
\]
This argument demonstrates the equivalence of Baxendale's setting (in terms
of $\gamma ^{\prime }(0)$ and $\gamma ^{\prime \prime }(0)$) and our setting
(in terms of $\lambda ^{\prime }(0)$ and $\lambda ^{\prime \prime }(0)$).

\section{First Derivative}

Let $y(x,\mu )$ denote the solution of the Ince BVP\ with $\lambda =\lambda
(\mu )$ that satisfies the initial conditions $y(0)=1$, $y^{\prime }(0)=0$.
Recall that $\lambda (\mu )$ is the leftmost such eigenvalue. Define 
\[
z(x,\mu )=\frac \partial {\partial \mu }y(x,\mu ). 
\]
Differentiate Ince's equation with respect to $\mu $, yielding 
\[
z^{\prime \prime }(x,\mu )+c\sin (2x)z^{\prime }(x,\mu )+\left( \lambda
^{\prime }(\mu )-c\cos (2x)\right) y(x,\mu )+\left( \lambda (\mu )-\mu
\,c\cos (2x)\right) z(x,\mu )=0. 
\]
Call this $(*)$. Set $\mu =0$, yielding 
\[
z^{\prime \prime }(x,0)+c\sin (2x)z^{\prime }(x,0)+\lambda ^{\prime
}(0)-c\cos (2x)=0 
\]
because $\lambda (0)=0$ and $y(x,0)=1$. Multiply both sides of $(*)$ by $%
\exp (-\frac c2\cos (2x))$, yielding 
\[
\left[ \exp (-\tfrac c2\cos (2x))z^{\prime }(x,0)\right] ^{\prime }=\left(
c\cos (2x)-\lambda ^{\prime }(0)\right) \exp (-\tfrac c2\cos (2x)). 
\]
Since $z^{\prime }(0,0)=0$, 
\[
\exp (-\tfrac c2\cos (2x))z^{\prime }(x,0)=\dint\limits_0^x\left( c\cos
(2\theta )-\lambda ^{\prime }(0)\right) \exp (-\tfrac c2\cos (2\theta
))d\theta . 
\]
Now $z^{\prime }(x,0)$ has period $\pi $, hence $z^{\prime }(\pi ,0)=0$ and
therefore 
\[
\dint\limits_0^\pi \left( c\cos (2\theta )-\lambda ^{\prime }(0)\right) \exp
(-\tfrac c2\cos (2\theta ))d\theta =0. 
\]
It follows that 
\[
\lambda ^{\prime }(0)=\frac{c\dint\limits_0^\pi \cos (2\theta )\exp (-\tfrac
c2\cos (2\theta ))d\theta }{\dint\limits_0^\pi \exp (-\tfrac c2\cos (2\theta
))d\theta }=c\frac{I_1(-\frac c2)}{I_0(-\frac c2)} 
\]
where $I_0$, $I_1$ are modified Bessel functions \cite{AS}. This reproduces
Baxendale's formula for $\ell $.

\section{Eigenfunction}

In the next section, we will need a formula for $z(x,0)$. Note that 
\[
y(x,\mu )\approx 1+\mu \,z(x,0) 
\]
when $\mu \approx 0$, thus a byproduct of our work is an approximation of
the leftmost eigenfunction of Ince's equation to first order.

From \cite{AS}, the following expressions hold: 
\[
\exp (-\tfrac c2\cos (2\theta ))=I_0(-\tfrac c2)+2\dsum\limits_{k=1}^\infty
I_k(-\tfrac c2)\cos (2k\,\theta ), 
\]
\[
\exp (-\tfrac c2\cos (2\theta ))\cos (2\theta )=I_1(-\tfrac
c2)+\dsum\limits_{k=1}^\infty \left( I_{k-1}(-\tfrac c2)+I_{k+1}(-\tfrac
c2)\right) \cos (2k\,\theta ) 
\]
hence the integrand $\left( c\cos (2\theta )-\lambda ^{\prime }(0)\right)
\exp (-\tfrac c2\cos (2\theta ))$ is equal to 
\begin{eqnarray*}
&&c\dsum\limits_{k=1}^\infty \left( I_{k-1}(-\tfrac c2)+I_{k+1}(-\tfrac
c2)\right) \cos (2k\,\theta )-2c\frac{I_1(-\tfrac c2)}{I_0(-\tfrac c2)}%
\dsum\limits_{k=1}^\infty I_k(-\tfrac c2)\cos (2k\,\theta ) \\
&=&c\dsum\limits_{k=1}^\infty \stackunder{j_k}{\underbrace{\left(
I_{k-1}(-\tfrac c2)-2\frac{I_1(-\tfrac c2)}{I_0(-\tfrac c2)}I_k(-\tfrac
c2)+I_{k+1}(-\tfrac c2)\right) }}\cos (2k\,\theta ).
\end{eqnarray*}
Integrating once gives 
\[
\exp (-\tfrac c2\cos (2x))z^{\prime }(x,0)=\frac c2\dsum\limits_{k=1}^\infty 
\frac{j_k}k\sin (2k\,x). 
\]
Integrating twice gives 
\[
z(x,0)=\frac c2\dsum\limits_{k=1}^\infty \frac{j_k}k\dint\limits_0^x\exp
(\tfrac c2\cos (2\theta ))\sin (2k\,\theta )d\theta . 
\]
Although we do not explicitly write out the integrals here, they are
elementary and can be computed symbolically for arbitrary $k$.

\section{Second Derivative}

Define 
\[
w(x,\mu )=\frac{\partial ^2}{\partial ^2\mu }y(x,\mu )=\frac \partial
{\partial \mu }z(x,\mu ). 
\]
Differentiate equation $(*)$ with respect to $\mu $, yielding 
\[
\begin{array}{c}
w^{\prime \prime }(x,\mu )+c\sin (2x)w^{\prime }(x,\mu )+\lambda ^{\prime
\prime }(\mu )y(x,\mu )+2\left( \lambda ^{\prime }(\mu )-c\cos (2x)\right)
z(x,\mu ) \\ 
+\left( \lambda (\mu )-\mu \,c\cos (2x)\right) w(x,\mu )=0.
\end{array}
\]
Set $\mu =0$, yielding 
\[
w^{\prime \prime }(x,0)+c\sin (2x)w^{\prime }(x,0)+\lambda ^{\prime \prime
}(0)+2\left( \lambda ^{\prime }(0)-c\cos (2x)\right) z(x,0)=0 
\]
because $\lambda (0)=0$ and $y(x,0)=1$. We again multiply both sides by $%
\exp (-\frac c2\cos (2x))$: 
\[
\left[ \exp (-\tfrac c2\cos (2x))w^{\prime }(x,0)\right] ^{\prime }=\left[
2\left( c\cos (2x)-\lambda ^{\prime }(0)\right) z(x,0)-\lambda ^{\prime
\prime }(0)\right] \exp (-\tfrac c2\cos (2x)). 
\]
Since $w^{\prime }(0,0)=0$, 
\[
\exp (-\tfrac c2\cos (2x))w^{\prime }(x,0)=\dint\limits_0^x\left[ 2\left(
c\cos (2\theta )-\lambda ^{\prime }(0)\right) z(\theta ,0)-\lambda ^{\prime
\prime }(0)\right] \exp (-\tfrac c2\cos (2\theta ))d\theta . 
\]
Now $w^{\prime }(x,0)$ has period $\pi $, hence $w^{\prime }(\pi ,0)=0$ and
therefore 
\[
\dint\limits_0^\pi \left[ 2\left( c\cos (2\theta )-\lambda ^{\prime
}(0)\right) z(\theta ,0)-\lambda ^{\prime \prime }(0)\right] \exp (-\tfrac
c2\cos (2\theta ))d\theta =0. 
\]
It follows that 
\[
\lambda ^{\prime \prime }(0)=\frac 2{\pi \,I_0(-\frac c2)}\dint\limits_0^\pi
\left( c\cos (2\theta )-\lambda ^{\prime }(0)\right) z(\theta ,0)\exp
(-\tfrac c2\cos (2\theta ))d\theta 
\]
where $z(\theta ,0)$ is defined via an infinite series in the preceding
section.

\section{Numerical Verification}

One way to compute $\lambda (\mu )$ is to construct the infinite tridiagonal
matrix \cite{Vm} 
\[
M=\left( 
\begin{array}{cccccc}
r_0 & \sqrt{2}q_{-1} & 0 & 0 & 0 &  \\ 
\sqrt{2}q_0 & r_1 & q_{-2} & 0 & 0 &  \\ 
0 & q_1 & r_2 & q_{-3} & 0 &  \\ 
0 & 0 & q_2 & r_3 & q_{-4} &  \\ 
0 & 0 &  & q_3 & r_4 & \ddots \\ 
&  &  &  & \ddots & \ddots
\end{array}
\right) 
\]
where $r_n=4n^2$ and $q_n=(-n+\mu /2)c$. The leftmost eigenvalue of $M$ is $%
\lambda (\mu )$.

Another way to compute $\lambda (\mu )$ is to solve the continued fraction
equation \cite{Vm} 
\[
-\frac \lambda 2=\dfrac{\left. p_0\;\;\;\;\;\;\;\;\;\;\right| }{\left|
4\cdot 1^2-\lambda \right. }-\dfrac{\left. p_1\;\;\;\;\;\;\;\;\;\;\right| }{%
\left| 4\cdot 2^2-\lambda \right. }-\dfrac{\left.
p_2\;\;\;\;\;\;\;\;\;\;\right| }{\left| 4\cdot 3^2-\lambda \right. }-\dfrac{%
\left. p_3\;\;\;\;\;\;\;\;\;\;\right| }{\left| 4\cdot 4^2-\lambda \right. }-%
\dfrac{\left. p_4\;\;\;\;\;\;\;\;\;\;\right| }{\left| 4\cdot 5^2-\lambda
\right. }-\cdots 
\]
where $p_n=(-n+\mu /2)(n+1+\mu /2)c^2$.

High precision estimates of the derivatives of $\lambda (\mu )$ at zero are
found via 
\[
\begin{array}{ccc}
\lambda ^{\prime }(0)\approx \dfrac{\lambda (\mu )}\mu , &  & \lambda
^{\prime \prime }(0)\approx \dfrac{\dfrac{\lambda (\mu )}\mu -c\dfrac{%
I_1(-\frac c2)}{I_0(-\frac c2)}}{\dfrac \mu 2}
\end{array}
\]
for $\mu \approx 0$. For example, when $a=1$, $b=-2$ and $\sigma =10$, we
obtain 
\[
\ell =-0.4887503163943852244580286...,
\]
\[
s^2=0.0112485762885419873084837...
\]
as the CLT\ parameter values. As another example, 
\[
\ell =0.3941998582469360577816389...,
\]
\[
s^2=0.3841476218435126147382099...
\]
when $\sigma =1$ (and $a$, $b$ remain unchanged). The same results were
obtained using exact formulation in sections 2 and 4, confirming our work.

\section{Acknowledgements}

Peter Baxendale suggested that I work on Ince's equation, Hans Volkmer
sketched out the derivation of $\lambda ^{\prime }(0)$ and $\lambda ^{\prime
\prime }(0)$, and Joe Keane found the generating function in \cite{AS} for $%
\exp (r\cos (\theta ))$. My grateful thanks go to all three!

\end{document}